\documentclass[12pt,reqno]{amsproc}
\usepackage[papersize={8.5in,11in},margin=1in]{geometry}

\usepackage{natbib,newtx,bm}
\usepackage[breaklinks=true,colorlinks,urlcolor=blue,citecolor=blue,linkcolor=blue]{hyperref}

\newtheorem{thm}{Theorem}[section]
\newtheorem{cor}[thm]{Corollary}
\newtheorem{lem}[thm]{Lemma}
\newtheorem{prop}[thm]{Proposition}
\theoremstyle{definition}

\theoremstyle{remark}
\newtheorem{rem}[thm]{\it Remark}
\newtheorem*{pf}{\it Proof}

\numberwithin{equation}{section}
\allowdisplaybreaks

\DeclareSymbolFont{CM}{OMX}{cmex}{m}{n}
\DeclareMathSymbol{\sumop}{\mathop}{CM}{"50}
\renewcommand{\sum}{\sumop}

\AtBeginDocument{.}

\begin{document}

\title[Multiple solutions for $(p,q)$-Laplacian equations]{\boldmath On $p(x)$-Laplacian equations in $\mathbb{R}^{N}$ with nonlinearity\\ sublinear
at zero}
\author[S. Liu \& C. Zhao]{Shibo Liu$^{\rm(a)}$ \& Chunshan Zhao$^{\rm(b)}$\vspace{-1em}}
\dedicatory{$^{\rm(a)}$Department of Mathematics \& Systems Engineering, Florida Institute of Technology\\Melbourne, FL 32901, USA\\[1ex]
$^{\rm(b)}$Department of Mathematical Sciences, Georgia Southern University\\Statesboro, GA 30460, USA}
\thanks{\emph{Email addresses}: \texttt{sliu@fit.edu} (S. Liu), \texttt{czhao@georgiasouthern.edu} (C. Zhao)}


\begin{abstract}
Let $p,q$ be functions on $\mathbb{R}^{N}$ satisfying $1\ll q\ll p\ll N$, we consider $p(x)$-Laplacian problems of the form
\[
\left\{
\begin{array}
[c]{l}%
-\Delta_{p(x)}u+V(x)\vert u\vert ^{p(x)-2}u=\lambda\vert
u\vert ^{q(x)-2}u+g(x,u)\text{,}\\
u\in W^{1,p(x)}(\mathbb{R}^{N})\text{.}%
\end{array}
\right.  
\]
To apply variational methods, we introduce a subspace  $X$ of $W^{1,p(x)}(\mathbb{R}^N)$ as our working space. Compact embedding from $X$ into $L^{q(x)}(\mathbb{R}^N)$ is proved, this enable us to get nontrivial solution of the problem; and two sequences of solutions going to $\infty$ and $0$ respectively, when $g(x,\cdot)$ is odd.

\end{abstract}
\maketitle

\section{Introduction}

For functions $a,b:\mathbb{R}^{N}\rightarrow\mathbb{R}$ we write $a\ll b$ if
$\inf_{\mathbb{R}^{N}}(  a-b)  <0$,%
\[
a_{-}=\inf_{\mathbb{R}^{N}}a\text{,\qquad}a_{+}=\sup_{\mathbb{R}^{N}}a\text{.}%
\]
Let $p\in C^{0,1}(\mathbb{R}^{N})$, $q\in L^\infty(\mathbb{R}^{N})$ be such that $1\ll q\ll p\ll N$, we
consider the following problem%
\begin{equation}
\left\{
\begin{array}
[c]{l}%
-\Delta_{p(x)}u+V(x)\vert u\vert ^{p(x)-2}u=\lambda\vert
u\vert ^{q(x)-2}u+g(x,u)\text{,}\\
u\in W^{1,p(x)}(\mathbb{R}^{N})\text{,}%
\end{array}
\right.  \label{e2}%
\end{equation}
where $\Delta_{p(x)}u=\operatorname{div}(\vert \nabla u\vert
^{p(x)-2}\nabla u)$ is the $p(x)$-Laplacian of $u$, $W^{1,p(x)}(\mathbb{R}^{N})$ is the variable exponent
Sobolev space that will be recalled in the next section.

By
definition, $u\in W^{1,p(x)}(\mathbb{R}^{N})$ is a (weak) solution of
(\ref{e2}) if%
\[
\int\left(  \left\vert \nabla u\right\vert ^{p(x)-2}\nabla u\cdot\nabla
h+V(x)\vert u\vert ^{p(x)-2}uh\right)  =\lambda\int\left\vert
u\right\vert ^{q(x)-2}uh+\int g(x,u)h
\]
for all test functions $h\in C_{0}^{\infty}(\mathbb{R}^{N})$, where from now on
all integrals are taken over $\mathbb{R}^{N}$ unless stated explicitly.

The variable exponent problems like (\ref{e2}) arise in some applications such as elasticity, electrorheological fluids, flow in porous media, see
\cite{MR1329546,MR1810360,MR2790542}. 
Such problems have received considerable attention in the last two decades, see e.g.\  \cite{MR2394068,MR2674092,MR3419907,MR4218373}. In these papers $\lambda=0$ and $g(x,t)$ behaves like $o(|t|^{p(x)-1})$ as $t\to0$, which imply that $u=0$ is a local minimizer of the variational functional $\Phi_\lambda$ given in  (\ref{O}), and the Mountain Pass Theorem of \cite{MR0370183} can be applies.

Since $q\ll p$, if $\lambda>0$, that is, the problem (\ref{e2}) is $p(x)$-sublinear at zero, then $u=0$ is not a local minimizer of $\Phi_\lambda$ any more. There are also some works on $p(x)$-Laplacian equation (\ref{e2}) on $\mathbb{R}^N$ whose nonlinearity (the right hand side of the equation) is $p(x)$-sublinear, see e.g.\ \cite{MR3725473}, \citet[Theorem 1.2]{MR3884271} and \citet[Theorems 3.1 \& 3.2]{MR2092084}.

For $(p,q)\in C^{0,1}\times L^\infty$ satisfying $1\ll q\ll p\ll N$, to study the problem (\ref{e2}) we assume:

\begin{enumerate}
\item[$\left(  V\right)  $] $V\in C(\mathbb{R}^{N})$, $V_{-}>0$,
$V^{-q(\cdot)/\left(  p(\cdot)-q(\cdot)\right)  }\in L^{1}(\mathbb{R}^{N})$.

\item[$\left(  g_{0}\right)  $] $g\in C(\mathbb{R}^{N}\times\mathbb{R})$
satisfies%
\begin{equation}
\left\vert g(x,t)\right\vert \leq C\big(  \vert t\vert
^{p(x)-1}+\left\vert t\right\vert ^{s(x)-1}\big)  \label{e5}%
\end{equation}
for some $s\in L^\infty(\mathbb{R}^{N})$ with $p\le s\ll p^{\ast}$, where $p^{\ast
}\in L^\infty(\mathbb{R}^{N})$ is given in (\ref{e3}).

\item[$\left(  g_{1}\right)  $] There is a constant $\mu>p_{+}$ such that%
\[
0<G(x,t):=\int_{0}^{t}g(x,\cdot)\leq\frac{1}{\mu}tg(x,t)\qquad\text{for
$x\in\mathbb{R}^{N}$, $t\neq0$.}%
\]

\end{enumerate}
The condition $(g_1)$ is the well-known Ambrosetti-Rabinowitz (AR) condition introduced in \cite{MR0370183} for the case $p(x)\equiv2$. Due to the present of the $p(x)$-sublinear term $\lambda|t|^{q(x)-2}t$, our nonlinearity (the right hand side of (\ref{e2})) does not satisfy (AR).
\begin{thm}
\label{t1}Suppose $\left(  V\right)  $, $\left(  g_{0}\right)  $ and $\left(
g_{1}\right)  $ are satisfied, $p_{-}>q_{+}$.

\begin{enumerate}
\item If in addition%
\begin{equation}
\lim_{\left\vert t\right\vert \rightarrow0}\frac{G(x,t)}{\left\vert
t\right\vert ^{p(x)}}=0\text{\qquad uniformly in }x\in\mathbb{R}^{N}\text{,}
\label{e4}%
\end{equation}
then for some $\lambda_{0}>0$, the problem \eqref{e2} has a nontrivial
solution provided $\lambda\in\left(  -\lambda_{0},\lambda_{0}\right)  $.

\item If in addition $g(x,\cdot)$ is odd for\ $x\in\mathbb{R}^{N}$, then for
all $\lambda\in\mathbb{R}$ the problem \eqref{e2} has a sequence of solutions
$\left\{  u_{n}\right\}  $ such that%
\begin{equation}
\int\frac{1}{p(x)}\left(  \left\vert \nabla u_{n}\right\vert ^{p(x)}%
+V(x)\vert u_{n}\vert ^{p(x)}\right)  -\int\frac{\lambda\vert u_{n}\vert ^{q(x)}+ G(x,u_{n})}%
{q(x)}\rightarrow
+\infty\text{.} \label{xd}%
\end{equation}

\end{enumerate}
\end{thm}

\begin{rem}
If $p(x)\equiv2$, the $p(x)$-Laplacian $\Delta_{p(x)}$ reduces to the classical Laplacian and the left hand side of (\ref{e2}) becomes the Schr\"{o}dinger operator $S=-\Delta+V$. The semilinear nature of this situation allows us to deal with the indefinite case that $V_-<0$, see \cite{doi:10.1080/17476933.2023.2253533}.
\end{rem}

We point out that if $\lambda>0$, then without assuming $\left(  g_{1}\right)
$, which is crucial for getting the boundedness of Palais-Smale sequences (see the proof of Lemma \ref{l2}),
the problem (\ref{e2}) has another sequence of solutions $\left\{  u_{n}\right\}  $
such that the left hand side of (\ref{xd}) is negative and converging to zero.
In fact, we have a more general result. Consider the following problem%
\begin{equation}
\left\{
\begin{array}
[c]{l}%
-\Delta_{p(x)}u+V(x)\vert u\vert ^{p(x)-2}u=f(x,u)\text{,}\\
u\in W^{1,p(x)}(\mathbb{R}^{N})\text{,}%
\end{array}
\right.  \label{2e}%
\end{equation}
with the following assumption on $f\in C(\mathbb{R}^{N}\times\mathbb{R})$:

\begin{enumerate}
\item[$\left(  f_{0}\right)  $] $\left\vert f(x,t)\right\vert \leq C\big(\vert t\vert
^{q(x)-1}+\vert t\vert ^{s(x)-1}\big)$ for some $s\in L^\infty(\mathbb{R}^{N})$ with $p\le s\ll p^{\ast}$, and%
\begin{equation}
\lim_{\left\vert t\right\vert \rightarrow0}\frac{F(x,t)}{\left\vert
t\right\vert ^{p_{-}}}=+\infty\text{\quad uniformly on }\mathbb{R}%
^{N}\text{,\quad where }F(x,t)=\int_{0}^{t}f(x,\cdot)\text{.} \label{f2}%
\end{equation}

\end{enumerate}

\begin{thm}
\label{t2}Suppose $\left(  V\right)  $ and $\left(  f_{0}\right)  $ are
satisfied. If in addition $f(x,\cdot)$ is odd for\ $x\in\mathbb{R}^{N}$, then
the problem \eqref{2e} has a sequence of solutions $\left\{  u_{n}\right\}  $
such that $\Phi(u_{n})\leq0$ and $\Phi(u_{n})\rightarrow0$, where%
\[
\Phi(u)=\int\frac{1}{p(x)}\left(  \left\vert \nabla u\right\vert
^{p(x)}+V(x)\vert u\vert ^{p(x)}\right)  -\int F(x,u)\text{.}%
\]
\end{thm}

The paper is organized as follow. In Section 2 we briefly recall some notions and results on variable exponent Lebesgue space $L^{p(x)}(\mathbb{R}^N)$ and Sobolev space $W^{1,p(x)}(\mathbb{R}^N)$. A suitable subspace $X$ of $W^{1,p(x)}(\mathbb{R}^N)$ is introduced and we show that there is a compact embedding $X\hookrightarrow L^{q(x)}(\mathbb{R}^N)$. The proofs of Theorems \ref{t1} and \ref{t2} are presented in Section 3 and Section 4, respectively.

\section{Variable exponent Lebesgue and Sobolev spaces}

To study the problem (\ref{e2}) we recall the variable exponent Lebesgue
spaces and Sobolev spaces, see \cite{MR1866056,MR1859337} for more details.

Let $\Omega$ be an open smooth subset of $\mathbb{R}^{N}$, $p\in
L^\infty(\Omega)$ be such that $1\le p\ll N$. On the vector space%
\[
L^{p(x)}(\Omega)=\left\{  u:\Omega\rightarrow\mathbb{R}\left\vert \int%
_{\Omega}\left\vert u\right\vert ^{p(x)}<\infty\right.  \right\}
\]
we equip the Luxembury norm%
\[
\left\vert u\right\vert _{p(x);\Omega}=\inf\left\{  \xi>0\left\vert
\int_{\Omega}\left\vert \frac{u}{\xi}\right\vert ^{p(x)}\leq1\right.
\right\}  \text{.}%
\]
Then the variable exponent Lebesgue space $L^{p(x)}(\Omega)$ becomes a
separable uniformly convex Banach space. If $\Omega=\mathbb{R}^{N}$ we simply
write $\left\vert u\right\vert _{p(x)}$ for $\left\vert u\right\vert
_{p(x);\mathbb{R}^{N}}$. For $u\in L^{p(x)}(\Omega)$ we have the following
inequality%
\begin{equation}
\left\vert u\right\vert _{p(x);\Omega}\leq\left(  \int_{\Omega}\left\vert
u\right\vert ^{p(x)}\right)  ^{1/p_{+}}+\left(  \int_{\Omega}\left\vert
u\right\vert ^{p(x)}\right)  ^{1/p-}\text{.} \label{w}%
\end{equation}

The H\"{o}lder conjugate function $\tilde{p}$ and the critical Sobolev
exponent function $p^{\ast}$ are defined via%
\begin{equation}
\tilde{p}(x)=\frac{p(x)}{p(x)-1}\text{,\qquad}p^{\ast}(x)=\frac
{Np(x)}{N-p(x)}\text{.} \label{e3}%
\end{equation}
Here, to avoid ambiguous notations $p'_\pm$, instead of $p'$ we use $\tilde{p}$ for H\"{o}lder conjugate. We have the H\"{o}lder inequality: for $u\in L^{p(x)}(\Omega)$ and $v\in
L^{\tilde{p}(x)}(\Omega)$ there holds
\begin{equation}
\int_{\Omega}\left\vert uv\right\vert \leq2\left\vert u\right\vert
_{p(x);\Omega}\left\vert v\right\vert _{\tilde{p}(x);\Omega}\text{.}
\label{eh}%
\end{equation}
The variable exponent Sobolev space $W^{1,p(x)}(\mathbb{R}^{N})$ is the
completion of $C_{0}^{\infty}(\mathbb{R}^{N})$ with respect to the norm%
\[
\left\Vert u\right\Vert _{1}=\left\vert u\right\vert _{p(x)}+\left\vert \nabla
u\right\vert _{p(x)}\text{.}%
\]

\begin{prop}[{\citet[Theorems 1.1 \& 1.3]{MR1859337}}]\label{ip}
Let $p\in C^{0,1}(\mathbb{R}^{N})$, $s\in L^\infty(\mathbb{R}^{N})$.

\begin{enumerate}
\item If $p\leq s\leq p^{\ast}$, then we have a continuous embedding $W^{1,p(x)}%
(\mathbb{R}^{N})\hookrightarrow L^{s(x)}(\mathbb{R}^{N})$.

\item \label{it2}If $1\leq s\ll p^{\ast}$, then we have a compact embedding
$W^{1,p(x)}(\mathbb{R}^{N})\hookrightarrow L_{\mathrm{loc}}^{s(x)}%
(\mathbb{R}^{N})$.
\end{enumerate}
\end{prop}

As a consequence of (2), if $u_{n}\rightharpoonup0$ in $W^{1,p(x)}(\mathbb{R}^{N})$, then for any given $R>0$ we have $|u_n|_{s(x);B_R}\to0$, which is equivalent to
\[
\int_{B_R}|u_n|^{s(x)}\to0\text{,}
\]
where $B_{R}$ is the $R$-ball in proper space (here is $\mathbb{R}^{N}$). We also denote by $B_{R}^{\mathrm{c}}$ its complement.

Let $X$ be the completion of $C_{0}^{\infty}(\mathbb{R}^{N})$ with respect to
the Luxemburg norm%
\[
\left\Vert u\right\Vert =\inf\left\{  \xi>0\left\vert \int\left(  \left\vert
\frac{\nabla u}{\xi}\right\vert ^{p(x)}+V(x)\left\vert \frac{u}{\xi
}\right\vert ^{p(x)}\right)  \leq1\right.  \right\}  \text{.}%
\]
Then $X$ is a separable uniformly convex Banach space which can be continuously embedded into
$W^{1,p(x)}(\mathbb{R}^{N})$ due to the assumption $V_{-}>0$.

The following lemma is generalization of \citet[Lemma 2.1]{doi:10.1080/17476933.2023.2253533}, which deals with the semilinear case that $p(x)\equiv2$ and $q(x)$ is a constant in $(1,2)$.

\begin{lem}
\label{l1} Assuming $\left(  V\right)  $, then there is a continuous and compact
embedding $X\hookrightarrow L^{q(x)}(\mathbb{R}^{N})$.
\end{lem}

\begin{pf}
Given $\varepsilon\in\left(  0,1\right)  $, because $V^{-q(\cdot)/\left(
p(\cdot)-q(\cdot)\right)  }\in L^{1}(\mathbb{R}^{N})$ and%
\[
\frac{p(x)}{p(x)-q(x)}\leq\frac{p_{+}}{\inf_{\mathbb{R}^{N}}\left(
p-q\right)  }=:\alpha
\]
for all $x\in\mathbb{R}^{N}$, there is $R>0$ such that%
\begin{equation}
\int_{\left\vert x\right\vert \geq R}\frac{V^{-q(x)/\left(  p(x)-q(x)\right)
}(x)}{\varepsilon^{p(x)/\left(  p(x)-q(x)\right)  }}\leq\frac{1}%
{\varepsilon^{\alpha}}\int_{\left\vert x\right\vert \geq R}V^{-q(x)/\left(
p(x)-q(x)\right)  }(x)\leq1\text{.} \label{ep}%
\end{equation}
We deduce%
\begin{align}
&\hspace{-3em}\left\vert V^{-q(x)/p(x)}\right\vert _{p(x)/\left(  p(x)-q(x)\right)
;B_{R}^{\mathrm{c}}}  \nonumber\\&  =\inf\left\{  \xi>0\left\vert \int_{\left\vert
x\right\vert \geq R}\left\vert \frac{V^{-q(x)/p(x)}(x)}{\xi}\right\vert
^{p(x)/\left(  p(x)-q(x)\right)  }\leq1\right.  \right\} \nonumber\\
&  =\inf\left\{  \xi>0\left\vert \int_{\left\vert x\right\vert \geq R}%
\frac{V^{-q(x)/\left(  p(x)-q(x)\right)  }(x)}{\xi^{p(x)/\left(
p(x)-q(x)\right)  }}\leq1\right.  \right\}  \leq\varepsilon\label{e}%
\end{align}
because by (\ref{ep}), $\varepsilon$ belongs to the set over which we take infimum.

Let $u\in X$. We have%
\begin{align*}
\int_{B_{R}^{\mathrm{c}}}V(x)\vert u\vert ^{p(x)}  &  \leq
\int\left(  \left\vert \nabla u\right\vert ^{p(x)}+V(x)\vert u\vert
^{p(x)}\right) \\
&  \leq1+\left\Vert u\right\Vert ^{p_{+}}+\left\Vert u\right\Vert ^{p_{-}%
}\text{.}%
\end{align*}
This and $(  q_{+}/p_{-})  (  p/q)  \geq1$ yield%
\[
\int_{B_{R}^{\mathrm{c}}}\left\vert \frac{V^{q(x)/p(x)}(x)\vert
u\vert ^{q(x)}}{\left(  1+\left\Vert u\right\Vert ^{p_{+}}+\left\Vert
u\right\Vert ^{p_{-}}\right)  ^{q_{+}/p_{-}}}\right\vert ^{p(x)/q(x)}\leq
\int_{B_{R}^{\mathrm{c}}}\frac{V(x)\vert u\vert ^{p(x)}%
}{1+\left\Vert u\right\Vert ^{p_{+}}+\left\Vert u\right\Vert ^{p_{-}}}%
\leq1\text{.}%
\]
Consequently%
\begin{align}
&\hspace{-3em}\left\vert V^{q(x)/p(x)}\vert u\vert ^{q(x)}\right\vert
_{p(x)/q(x);B_{R}^{\mathrm{c}}} \nonumber\\
&  =\inf\left\{  \xi>0\left\vert \int%
_{B_{R}^{\mathrm{c}}}\left\vert \frac{V^{q(x)/p(x)}(x)\vert u\vert
^{q(x)}}{\xi}\right\vert ^{p(x)/q(x)}\leq1\right.  \right\} \nonumber\\
&  \leq\big(  1+\left\Vert u\right\Vert ^{p_{+}}+\left\Vert u\right\Vert
^{p_{-}}\big)  ^{q_{+}/p_{-}}\text{.} \label{e1}%
\end{align}

It follows from (\ref{e}), (\ref{e1}) and the H\"{o}lder inequality (\ref{eh})
that
\begin{align}
\int_{\left\vert x\right\vert \geq R}\left\vert u\right\vert ^{q(x)}  &
=\int_{\left\vert x\right\vert \geq R}V^{q(x)/p(x)}(x)\vert u\vert
^{q(x)}\cdot V^{-q(x)/p(x)}(x)\nonumber\\
&  \leq2\left\vert V^{q(x)/p(x)}\vert u\vert ^{q(x)}\right\vert
_{p(x)/q(x);B_{R}^{\mathrm{c}}}\left\vert V^{-q(x)/p(x)}\right\vert
_{p(x)/\left(  p(x)-q(x)\right)  ;B_{R}^{\mathrm{c}}}\nonumber\\
&  \leq2\varepsilon\cdot\big(  1+\left\Vert u\right\Vert ^{p_{+}%
}+\left\Vert u\right\Vert ^{p_{-}}\big)  ^{q_{+}/p_{-}}\text{.} \label{ex}%
\end{align}
Due to the continuous embeddings $X\hookrightarrow W^{1,p(x)}(\mathbb{R}%
^{N})\hookrightarrow L^{q(x)}(B_{R})$, where the second embedding is given by
$u\mapsto u|_{B_{R}}$, we see that $u|_{B_{R}}\in L^{q(x)}(B_{R})$. Combining
(\ref{ex}) we deduce
\[
\int\left\vert u\right\vert ^{q(x)}=\int_{\left\vert x\right\vert
<R}\left\vert u\right\vert ^{q(x)}+\int_{\left\vert x\right\vert \geq
R}\left\vert u\right\vert ^{q(x)}<\infty\text{,}%
\]
that is $u\in L^{q(x)}(\mathbb{R}^{N})$. Therefore, $X\subset L^{q(x)}%
(\mathbb{R}^{N})$.

Now, let $\left\{  u_{n}\right\}  $ be a sequence in $X$ such that
$u_{n}\rightharpoonup0$ in $X$. Given $\varepsilon\in\left(  0,1\right)  $ we
choose $R>0$ satisfying (\ref{e}). Because $X\hookrightarrow W^{1,p(x)}%
(\mathbb{R}^{N})$ is continuous and $W^{1,p(x)}(\mathbb{R}^{N})\hookrightarrow
L_{\mathrm{loc}}^{q(x)}(\mathbb{R}^{N})$ is compact (see Proposition \ref{ip}), we have%
\[
\int_{\left\vert x\right\vert <R}\left\vert u_{n}\right\vert ^{q(x)}%
\rightarrow0\text{.}%
\]
This and (\ref{ex}) with $u$ replaced by $u_{n}$ yield%
\begin{align*}
\varlimsup_{n\rightarrow\infty}\int\left\vert u_{n}\right\vert ^{q(x)}  &
=\varlimsup_{n\rightarrow\infty}\left(  \int_{\left\vert x\right\vert \geq
R}+\int_{\left\vert x\right\vert <R}\right)  \left\vert u_{n}\right\vert
^{q(x)}\\
&  \leq2\varepsilon\cdot\sup_{n}\left(  1+\left\Vert u_{n}\right\Vert ^{p_{+}%
}+\left\Vert u_{n}\right\Vert ^{p_{-}}\right)  ^{q_{+}/p_{-}}\text{.}%
\end{align*}
Letting $\varepsilon\rightarrow0$ we deduce (noting $\sup_{n}\left\Vert
u_{n}\right\Vert <\infty$)%
\[
\int\left\vert u_{n}\right\vert ^{q(x)}\rightarrow0\text{,}%
\]
which implies that $u_{n}\rightarrow0$ in $L^{q(x)}(\mathbb{R}^{N})$. This
proves that the embedding $X\hookrightarrow L^{q(x)}(\mathbb{R}^{N})$ is
continuous and compact.\qed
\end{pf}

\begin{cor}
\label{c1}Let $s\in L^\infty(\mathbb{R}^{N})$ be such that $q\ll s\ll p^{\ast}%
$. Then there is a continuous and compact embedding $X\hookrightarrow
L^{s(x)}(\mathbb{R}^{N})$.
\end{cor}

\begin{pf}
This follows from an interpolation argument as in \citet[Page 2570]{MR2674092}.
There is $\lambda:\mathbb{R}^{N}\rightarrow\left(  0,1\right)  $ such that%
\[
\frac{1}{s(x)}=\frac{\lambda(x)}{q(x)}+\frac{1-\lambda(x)}{p^{\ast}%
(x)}\text{\qquad a.e.\  $x\in\mathbb{R}^{N}$.}%
\]
For continuous $\alpha:\mathbb{R}^{N}\rightarrow\left(  1,\infty\right)  $ and
its H\"{o}lder conjugate $\tilde{\alpha}:\mathbb{R}^{N}\rightarrow\left(
0,\infty\right)  $ given by%
\[
\alpha(x)=\frac{q(x)}{s(x)\lambda(x)}\text{,\qquad}\tilde{\alpha}(x)=\frac{p^{\ast}%
(x)}{s(x)\left(  1-\lambda(x)\right)  }\text{,}%
\]
we deduce from (\ref{w}) and the H\"{o}lder inequality (\ref{eh})%
\begin{align*}
\int\left\vert u\right\vert ^{s(x)}  &  =\int\left\vert u\right\vert
^{q(x)/\alpha(x)}\left\vert u\right\vert ^{p^{\ast}(x)/\tilde{\alpha}(x)}%
\leq2\left\vert \left\vert u\right\vert ^{q(x)/\alpha(x)}\right\vert
_{\alpha(x)}\left\vert \left\vert u\right\vert ^{p^{\ast}(x)/\tilde{\alpha}
(x)}\right\vert _{\tilde{\alpha}(x)}\\
&  \leq2\left(  \left(  \int\left\vert u\right\vert ^{q(x)}\right)
^{1/\alpha_{+}}+\left(  \int\left\vert u\right\vert ^{q(x)}\right)
^{1/\alpha_{-}}\right) \\
&  \qquad\qquad\times\left(  \left(  \int\left\vert u\right\vert ^{p^{\ast
}(x)}\right)  ^{1/\tilde{\alpha}_{+}}+\left(  \int\left\vert u\right\vert ^{p_{\ast
}(x)}\right)  ^{1/\tilde{\alpha}_{-}}\right)  \text{.}%
\end{align*}
The desired result follows from this and Lemma \ref{l1}.\qed
\end{pf}

\begin{prop}
\label{p1}Let $f\in C(\mathbb{R}^{N}\times\mathbb{R})$ satisfy%
\begin{equation}
\left\vert f(x,t)\right\vert \leq C\big(\vert t\vert
^{q(x)-1}+\left\vert t\right\vert ^{s(x)-1}\big)  \text{,} \label{g}%
\end{equation}
where $s\in L^\infty(\mathbb{R}^{N})$ satisfies $q\ll s\ll p^{\ast}%
$. If $u_{n}\rightharpoonup u$ in $X$, then
\[
\int f(x,u_{n})(  u_{n}-u)  \rightarrow0\text{,}\qquad\int F(x,u_n)\to\int F(x,u)\text{.}
\]

\end{prop}

\begin{pf}
From Lemma \ref{l1} and Corollary \ref{c1},%
\[
\left\vert u_{n}-u\right\vert _{q(x)}\rightarrow0\text{,\qquad}\left\vert
u_{n}-u\right\vert _{s(x)}\rightarrow0\text{.}%
\]
Because $\left\{  u_{n}\right\}  $ is bounded in both $L^{q(x)}(\mathbb{R}%
^{N})$ and $L^{s(x)}(\mathbb{R}^{N})$, using (\ref{w}) we get%
\begin{align*}
\left\vert \left\vert u_{n}\right\vert ^{q(x)-1}\right\vert _{\tilde{q}(x)}
&  \leq\left(  \int\left\vert \left\vert u_{n}\right\vert ^{q(x)-1}\right\vert
^{\tilde{q}(x)}\right)  ^{1/\tilde{q}_{+}}+\left(  \int\left\vert \left\vert
u_{n}\right\vert ^{q(x)-1}\right\vert ^{\tilde{q}(x)}\right)  ^{1/\tilde{q}_{-}%
}\\
&  =\left(  \int\left\vert u_{n}\right\vert ^{q(x)}\right)  ^{1/\tilde{q}_{+}}+\left(  \int\left\vert u_{n}\right\vert ^{q(x)}\right)  ^{1/\tilde{q}_{-}%
}\leq B\text{,}%
\end{align*}
and%
\[
\left\vert \left\vert u_{n}\right\vert ^{s(x)-1}\right\vert _{\tilde{s}%
(x)}\leq\left(  \int\left\vert u_{n}\right\vert ^{s(x)}\right)  ^{1/\tilde{s}_{+}%
}+\left(  \int\left\vert u_{n}\right\vert ^{s(x)}\right)
^{1/\tilde{s}_{-}}\leq B
\]
for some $B>0$. By the growth condition (\ref{g}) and H\"{o}lder inequality
(\ref{eh}), we have%
\begin{align}
&\hspace{-3em}\left\vert \int f(x,u_{n})(  u_{n}-u)  \right\vert \nonumber\\
&  \leq
C\left(  \int\left\vert u_{n}\right\vert ^{q(x)-1}\left\vert u_{n}%
-u\right\vert +\int\left\vert u_{n}\right\vert ^{s(x)-1}\left\vert
u_{n}-u\right\vert \right) \nonumber\\
&  \leq2C\left(  \left\vert \left\vert u_{n}\right\vert ^{q(x)-1}\right\vert
_{\tilde{q}(x)}\left\vert u_{n}-u\right\vert _{q(x)}+\left\vert \left\vert
u_{n}\right\vert ^{s(x)-1}\right\vert _{\tilde{s}(x)}\left\vert
u_{n}-u\right\vert _{s(x)}\right) \nonumber\\
&  \leq2CB\left(  \left\vert u_{n}-u\right\vert _{q(x)}+\left\vert
u_{n}-u\right\vert _{s(x)}\right)  \rightarrow0\text{.}\label{po}%
\end{align}
To prove the second assertion, we apply the mean value theorem to $\Psi:u\mapsto\int F(x,u)$. There is $h_n\in [u,u_n]$, the segment in $X$ jointing $u$ and $u_n$, such that
\begin{align*}
\left|\int F(x,u_n)-\int F(x,u)\right|&=|\langle\Psi'(h_n),u_n-u\rangle|\\
&=\left|\int f(x,h_n)(u_n-u)\right|\text{.}
\end{align*}
Because $\{h_n\}$ is bounded, the right hand side goes to zero as in (\ref{po}).\qed
\end{pf}

\section{Proof of Theorem \ref{t1}}

We consider the energy functional $\Phi_{\lambda}:X\rightarrow\mathbb{R}$,%
\begin{align}
\Phi_{\lambda}(u)=\int\frac{1}{p(x)}\left(  \left\vert \nabla u\right\vert
^{p(x)}+V(x)\vert u\vert ^{p(x)}\right)  -\lambda\int\frac{1}%
{q(x)}\left\vert u\right\vert ^{q(x)}-\int G(x,u)\text{.}\label{O}%
\end{align}
It is standard to show that $\Phi_{\lambda}$ is of class $C^{1}$,%
\begin{align*}
\left\langle \Phi_{\lambda}(u),h\right\rangle &=\int\left(  \left\vert
\nabla u\right\vert ^{p(x)-2}\nabla u\cdot\nabla h+V(x)\vert
u\vert ^{p(x)-2}uh\right)\\
&  \hspace{3em}-\lambda\int\left\vert u\right\vert
^{q(x)-2}uh-\int g(x,u)h
\end{align*}
for $u,h\in X$. Therefore, critical points of $\Phi_{\lambda}$ are
solutions of problem (\ref{e2}).

\begin{lem}
\label{l2} Let $\left(  V\right)  $, $\left(  g_{0}\right)  $ and $\left(
g_{1}\right)  $ hold, $p_{-}>q_{+}$. Then for all $\lambda\in\mathbb{R}$, the
functional $\Phi_{\lambda}$ satisfies the Palais-Smale $\left(  PS\right)  $ condition.
\end{lem}

\begin{pf}
Denote by $\ell$ the norm of the embedding $X\hookrightarrow L^{q(x)}%
(\mathbb{R}^{N})$. From the definition of $\left\vert \cdot\right\vert
_{q(x)}$ we have%
\begin{equation}
\int\left\vert u\right\vert ^{q(x)}\leq\left(  1+\left\vert u\right\vert
_{q(x)}\right)  ^{q_{+}}\leq\left(  1+\ell\left\Vert u\right\Vert \right)
^{q_{+}}\text{\qquad for all }u\in X\text{.} \label{el}%
\end{equation}
Let $\left\{  u_{n}\right\}  \subset X$ be a $\left(  PS\right)  $ sequence of
$\Phi_{\lambda}$, that is $\Phi_{\lambda}^{\prime}(u_{n})\rightarrow0$,
$\Phi_{\lambda}(u_{n})\rightarrow c$ for some $c\in\mathbb{R}$. To show that
$\left\{  u_{n}\right\}  $ is bounded we may assume that $\left\Vert
u_{n}\right\Vert \geq1$, hence%
\begin{equation}
\int\left(  \left\vert \nabla u_{n}\right\vert ^{p(x)}+V(x)\vert
u_{n}\vert ^{p(x)}\right)  \geq\left\Vert u_{n}\right\Vert ^{p_{-}%
}\text{.} \label{ek}%
\end{equation}
Using condition $\left(  g_{1}\right)  $, (\ref{el}) and (\ref{ek}), for
$n\gg1$ we have%
\begin{align*}
c+1  &  \geq\Phi_{\lambda}(u_{n})-\frac{1}{\mu}\langle\Phi_{\lambda}^{\prime
}(u_{n}),u_{n}\rangle\\
&  =\int\left(  \frac{1}{p(x)}-\frac{1}{\mu}\right)  \left(  \left\vert \nabla
u_{n}\right\vert ^{p(x)}+V(x)\vert u_{n}\vert ^{p(x)}\right) \\
&  \qquad\qquad+\int\frac{u_{n}g(x,u_{n})-\mu G(x,u_{n})}{\mu}+\int\left(
\frac{1}{\mu}-\frac{\lambda}{q(x)}\right)  \left\vert u_{n}\right\vert
^{q(x)}\\
&  \geq\left(  \frac{1}{p_{+}}-\frac{1}{\mu}\right)  \int\left(  \left\vert
\nabla u_{n}\right\vert ^{p(x)}+V(x)\vert u_{n}\vert ^{p(x)}%
\right)  +\left(  \frac{1}{\mu}-\frac{\lambda}{q_{-}}\right)  \int\left\vert
u_{n}\right\vert ^{q(x)}\\
&  \geq\left(  \frac{1}{p_{+}}-\frac{1}{\mu}\right)  \left\Vert u_{n}%
\right\Vert ^{p_{-}}-\left\vert \frac{1}{\mu}-\frac{\lambda}{q_{-}}\right\vert
\left(  1+\ell\left\Vert u_{n}\right\Vert \right)  ^{q_{+}}\text{.}%
\end{align*}
Since $p_{-}>q_{+}$, it follows that $\left\{  u_{n}\right\}  $ is bounded in
$X$.

Up to a subsequence we may assume $u_{n}\rightharpoonup u$ in $X$. Using (\ref{e5}) we see that
\[
f(x,t):=\lambda\left\vert t\right\vert ^{q(x)-2}t+g(x,t)
\]
satisfies (\ref{g}), Proposition \ref{p1} yields%
\begin{equation}
\int f(x,u_{n})(  u_{n}-u)  \rightarrow0\text{.}\label{he}%
\end{equation}
Let $\mathcal{L}:X\rightarrow X^{\ast}$ be defined by.%
\begin{equation}
\left\langle \mathcal{L}u,\phi\right\rangle =\int\left(  \left\vert \nabla
u\right\vert ^{p(x)-2}\nabla u\cdot\nabla\phi+V(x)\vert u\vert
^{p(x)-2}u\phi\right)  \text{.}\label{lo}%
\end{equation}
Similar to \citet[Theorem 3.1 (ii)]{MR1954585} it can be shown that
$\mathcal{L}$ is $\left(  S_{+}\right)  $-type. That is,

\begin{itemize}
\item $\left\langle \mathcal{L}u_{n},u_{n}-u\right\rangle \rightarrow0$ and
$u_{n}\rightharpoonup u$ in $X$ imply $u_{n}\rightarrow u$ in $X$.
\end{itemize}
Using (\ref{he}) and $\Phi_{\lambda}^{\prime}(u_{n})\rightarrow0$, we deduce
\[
\left\langle \mathcal{L}u_{n},u_{n}-u\right\rangle =\langle\Phi_{\lambda
}^{\prime}(u_{n}),u_{n}-u\rangle+\int f(x,u_{n})(  u_{n}-u)
\rightarrow0\text{.}%
\]
It follows that $u_{n}\rightarrow u$ in $X$.\qed
\end{pf}

\subsection{Proof of Theorem \ref{t1} (1)}

It is standard to get%
\[
\int G(x,u)=o(\left\Vert u\right\Vert ^{p_{+}})\text{\qquad as }\left\Vert
u\right\Vert \rightarrow0
\]
from $\left(  g_{0}\right)  $ and (\ref{e4}). Since
\[
\int\frac{1}{p(x)}\left(  \left\vert \nabla u\right\vert ^{p(x)}%
+V(x)\vert u\vert ^{p(x)}\right)  \geq \frac{1}{p_+}\left\Vert u\right\Vert
^{p_{+}}%
\]
for $\left\Vert u\right\Vert \leq1$, we see that there is $r\in\left(
0,1\right)  $ such that%
\[
\eta:=\inf_{\left\Vert u\right\Vert =r}\Phi_{0}(u)>0\text{,}%
\]
where $\Phi_0$ is given in (\ref{O}) with $\lambda=0$.

Let $\lambda_{0}>0$ be such that%
\[
\kappa:=\eta-\lambda_{0}\left(  1+\ell r\right)  ^{q_{+}}>0\text{,}%
\]
where $\ell$ is the norm of the embedding $X\hookrightarrow L^{q(x)}%
(\mathbb{R}^{N})$. If $\left\vert \lambda\right\vert <\lambda_{0}$ and
$\left\Vert u\right\Vert =r$, using (\ref{el}) we have%
\begin{align*}
\Phi_{\lambda}(u)  &  =\Phi_{0}(u)-\lambda\int\left\vert u\right\vert
^{q(x)}\\
&  \geq\eta-\left\vert \lambda\right\vert \left(  1+\ell r\right)  ^{q_{+}%
}>\kappa\text{.}%
\end{align*}
Hence%
\begin{equation}
b=\inf_{\left\Vert u\right\Vert =r}\Phi_{\lambda}(u)>0\text{.} \label{mp}%
\end{equation}
Take $v\in X\backslash\left\{  0\right\}  $. Using $\left(  g_{1}\right)  $ it
is well known that $\Phi(tv)\rightarrow-\infty$ as $t\rightarrow\infty$. Since
$\Phi_{\lambda}$ satisfies $\left(  PS\right)  $ (see Lemma \ref{l2}), by the
Mountain Pass Theorem \cite{MR0370183} $\Phi_{\lambda}$ has a critical point
$u$ with $\Phi_{\lambda}(u)\geq b$. Hence $u\neq0$ and it is a nontrivial
solution of (\ref{e2}).

\subsection{Proof of Theorem \ref{t1} (2)}

To prove Theorem \ref{t1} (2), we need the following Symmetric Mountain Pass
Theorem of \cite{MR0370183}, see also \citet[Theorem
9.12]{MR845785}.

\begin{prop}
\label{p3}Let $X$ be an infinite dimensional Banach space, $\Phi\in C^{1}(X)$
is an even functional satisfying $\left(  PS\right)  $ condition, $\Phi(0)=0$.
If $X=Y\oplus Z$ with $\dim Y<\infty$, and $\Phi$ satisfies

\begin{enumerate}
\item For any finite dimensional subspace $W$ of $X$, there is $R_{W}%
\in\left(  0,\infty\right)  $ such that $\Phi\leq0$ on $W\backslash B_{R_{W}}$.

\item There are $\alpha,\delta\in\left(  0,\infty\right)  $ such that
$\Phi|_{Z\cap\partial B_{\delta}}\geq\alpha$,
\end{enumerate}
then $\Phi$ has a sequence of critical values $c_{n}\rightarrow\infty$.
\end{prop}

In the setting of Theorem \ref{t1} (2), our $\Phi_{\lambda}$ is even,
$\Phi_{\lambda}(0)=0$. From Lemma \ref{l2} $\Phi_{\lambda}$ satisfies $\left(
PS\right)  $. To apply Proposition \ref{p3} for getting a
sequence of critical points $\left\{  u_{n}\right\}  $ with $\Phi_{\lambda
}(u_{n})\rightarrow+\infty$, thus finishing the proof of Theorem \ref{t1} (2),
it suffices to show that:

\begin{enumerate}
\item[(a)] $\Phi_{\lambda}$ is anti-coercive on any finite dimensional
subspace $W$.

\item[(b)] There is a subspace $Z$ with $\operatorname*{codim}Z<\infty$ such
that%
\begin{equation}
\alpha:=\inf_{Z\cap\partial B_{1}}\Phi_{\lambda}>0\text{.} \label{ef}%
\end{equation}

\end{enumerate}

\paragraph{\it Verification of (a).}

We emphasize that the condition $(g_1)$ does \emph{not} imply
\[
G(x,t)\ge a_3|t|^\mu\qquad\text{for $|t|$ large. }
\]
A counterexample is $G(x,t)=a(x)|t|^{\mu}$ being $a:\mathbb{R}^N\to(0,\infty)$ continuous and decaying to zero at infinity. 
Therefore, we argue as in \citet[page 2577]{MR2674092}, where the condition $(g_1)$ is not assumed.

Let $\left\{  u_{n}\right\}  \subset W$ be
such that $\left\Vert u_{n}\right\Vert \rightarrow\infty$. Because $\dim W<\infty$, for some $v\in
W\backslash\left\{  0\right\}  $ we have $\left\Vert u_{n}\right\Vert
^{-1}u_{n}\rightarrow v$ in $W$. Using (\ref{el}), the
following consequence of $\left(  g_{1}\right)  $
\[
\lim_{\left\vert t\right\vert \rightarrow\infty}\frac{G(x,u_{n}(x))}%
{\left\vert u_{n}(x)\right\vert ^{p_{+}}}=+\infty
\]
for $x\in\left\{  v\neq0\right\}  $, and the Fatou's lemma we deduce%
\[
\Lambda_{n}:=\frac{1}{\left\Vert u_{n}\right\Vert ^{p_{+}}}\int\left(
\frac{\lambda\left\vert u_{n}\right\vert ^{q(x)}}{q(x)}+G(x,u_{n})\right)
\geq o_{n}(1)+\int_{v\ne0}\frac{G(x,u_{n})}{\left\Vert u_{n}\right\Vert
^{p_{+}}}\rightarrow+\infty\text{.}%
\]
Consequently%
\begin{align*}
\Phi_{\lambda}(u_{n})  &  \leq\frac{1}{p_{-}}\left\Vert u_{n}\right\Vert
^{p_{+}}-\int\left(  \frac{\lambda\left\vert u_{n}\right\vert ^{q(x)}}%
{q(x)}+G(x,u_{n})\right) \\
&  =\left\Vert u_{n}\right\Vert ^{p_{+}}\left(  \frac{1}{p_{-}}-\Lambda
_{n}\right)  \rightarrow-\infty
\end{align*}
as $n\to\infty$. Thus $\Phi_{\lambda}$ is anti-coercive on $W$.

\paragraph{\it Verification of (b).}

Let $X$ be a reflexive and separable Banach space. It is well known that there
exist sequences $\left\{  \phi_{n}\right\}  \subset X$ and $\left\{
f^{m}\right\}  \subset X^{\ast}$ such that $\left\langle f^{m},\phi
_{n}\right\rangle =\delta_{n}^{m}$, the Kronecker delta, and%
\[
X=\overline{\operatorname*{span}}\left\{  \phi_{n}\mid n\in\mathbb{N}\right\}
\text{,\qquad}X^{\ast}=\overline{\operatorname*{span}}^{\mathrm{w}^{\ast}%
}\{  f^{m}\mid m\in\mathbb{N}\}  \text{,}%
\]
see \cite{MR2359536}. For weakly continuous functional $\Psi:X\rightarrow\mathbb{R}$, set%
\[
\beta_{k}=\sup_{Z_{k}\cap\partial B_{1}}\left\vert \Psi\right\vert
\text{,\qquad}k\in\mathbb{N}\text{,}%
\]
where $Z_{k}=\overline{\operatorname*{span}}\left\{  \phi_{n}\mid n\geq
k\right\}  $. Then $\beta_{k}\rightarrow0$, see \citet[Lemma
3.3]{MR2092084}.

By Proposition \ref{p1} and (\ref{e5}), on our Sobolev space $X$, the functional $\Psi
:X\rightarrow\mathbb{R}$ defined via%
\[
\Psi(u)=\lambda
\int\frac{1}{q(x)}\left\vert u\right\vert ^{q(x)}+\int G(x,u)
\]
is weakly continuous. Using the result of Fan \& Han we just mentioned, there
is a sequence $\left\{  \beta_{k}\right\}  \subset\lbrack0,\infty)$, such that
$\beta_{k}\rightarrow0$ and%
\begin{equation}
\Psi(u)\leq\beta_{k}\text{,\qquad for all }u\in
Z_{k}\cap\partial B_{1}\text{.} \label{eg}%
\end{equation}
Because $\beta_{k}\rightarrow0$, for some $k_0\in\mathbb{N}$ there holds%
\[
\alpha:=\frac{1}{p_{+}}-\beta_{k_0}>0\text{.}%
\]
From this and (\ref{eg}), we deduce that for $u\in Z_{k_0}%
\cap\partial B_{1}$%
\begin{align*}
\Phi_{\lambda}(u)  &  =\int\frac{1}{p(x)}\left(  \left\vert \nabla
u\right\vert ^{p(x)}+V(x)\vert u\vert ^{p(x)}\right)  -\Psi(u)\\
&  \geq\frac{1}{p_{+}}\left\Vert u\right\Vert ^{p_{+}}-\beta_{k_0}=\frac{1}{p_{+}}- \beta_{k_0}=\alpha\text{,}%
\end{align*}
and (\ref{ef}) is verified with $Z=Z_{k_0}$.

\begin{rem}
To get infinitely many solutions for superlinear $p(x)$-Laplacian equations with odd nonlinearities, the Fountain Theorem of Bartsch \cite{MR1219237} is widely used, see for example \cite{MR3749216,MR3815570,MR4156782,MR2674092,MR1954585}.
It has also been applied in many other problems, such as Hamiltonian systems \cite{MR4061076,MR3564748} and Schr\"{o}dinger-Poisson systems \cite{MR2548724,MR3303004}. To apply the Fountain Theorem one needs to choose a sequence $\{\delta_k\}\subset(0,\infty)$ and verify
\[
\inf_{Z_k\cap \partial B_{\delta_k}}\Phi_\lambda\to+\infty\text{\qquad as $k\to\infty$,}
\]
which is much stronger than (\ref{ef}). Therefore, it seems more pleasurable to apply the Symmetric Mountain Pass
Theorem for such problems. On the other hand, we should point out that the Fountain Theorem can be applied to general $\mathcal{G}$-symmetric functionals (which reduce to even functionals when $\mathcal{G}=\mathbb{Z}_2$).
\end{rem}

\section{Proof of Theorem \ref{t2}}

Let $\Phi:X\rightarrow\mathbb{R}$ be the energy functional associated to problem (\ref{2e}),
defined by%
\[
\Phi(u)=\int\frac{1}{p(x)}\left(  \left\vert \nabla u\right\vert
^{p(x)}+V(x)\vert u\vert ^{p(x)}\right)  -\int F(x,u)\text{.}%
\]
To find a sequence of critical points of $\Phi$ as claimed in Theorem \ref{t2}, we adapt the approach
of \cite{liu2024multiple} (where a similar $p(x)$-Laplacian equation on a
bounded domain $\Omega$ is considered), which is an improvement of \cite{MR4162412} (see \citet[Remark 3.3]{liu2024multiple} for comparison), and is based on the following
version of the Clark's theorem.

\begin{prop}
[{\citet[Theorem 1.1]{MR3400440}}]\label{pp1}Let $X$ be a Banach space and
$\Psi\in C^{1}(X)$ be an even coercive functional satisfying the
$\left(  PS\right)  _{c}$ condition for $c\leq0$ and $\Psi(0)=0$. If for any
$k\in\mathbb{N}$ there is a $k$-dimensional subspace $W_{k}$ and $\delta
_{k}>0$ such that%
\begin{equation}
\sup_{W_{k}\cap \partial B_{\delta_{k}}}\Psi<0\text{,}\label{y}%
\end{equation}
then $\Psi$ has a sequence of critical points $u_{k}\neq0$ such that
$\Psi(u_{k})\leq0$, $u_{k}\rightarrow0$.
\end{prop}

As in \cite{MR4162412,liu2024multiple}, let $\phi:[0,\infty)\rightarrow\lbrack0,1]$ be a decreasing $C^{\infty}%
$-function such that $\left\vert \phi^{\prime}(t)\right\vert \leq2$,%
\[
\phi(t)=1\text{\quad for }t\in\left[  0,1\right]  \text{,\qquad}%
\phi(t)=0\text{\quad for }t\in\lbrack2,\infty)\text{.}%
\]
We consider the truncated functional $\Psi:X\rightarrow\mathbb{R}$,%
\[
\Psi(u)=\varrho(u)-\phi(\varrho(u))\int F(x,u)\text{,}%
\]
where $\varrho:X\rightarrow\mathbb{R}$ is defined by%
\[
\varrho(u)=\int\frac{1}{p(x)}\left(  \left\vert \nabla u\right\vert
^{p(x)}+V(x)\vert u\vert ^{p(x)}\right)  \text{.}%
\]
Note that $\varrho$ is $C^{1}$, and $\varrho^{\prime}(u)=\mathcal{L}u$ for the
$\mathcal{L}:X\rightarrow X^{\ast}$ given in (\ref{lo}).

The derivative of $\Psi$ is given by%
\begin{align}
\langle\Psi^{\prime}(u),h\rangle &  =\langle\varrho^{\prime}(u),h\rangle
-\phi(\varrho(u))\int f(x,u)h-\left(  \int F(x,u)\right)  \phi^{\prime}%
(\varrho(u))\langle\varrho^{\prime}(u),h\rangle\nonumber\\
&  =\left(  1-\left(  \int F(x,u)\right)  \phi^{\prime}(\varrho(u))\right)
\langle \varrho^{\prime}(u),h\rangle-\phi(\varrho(u))\int f(x,u)h\label{e0}%
\end{align}
for $u,h\in X$.

\begin{lem}
\label{l6}The functional $\Psi$ is coercive and satisfies $\left(  PS\right)
_{c}$ for $c\leq0$.
\end{lem}

\begin{pf}
The proof is similar to that of \citet[Lemma  3.2]{liu2024multiple}. For the reader's convenience we sketch it  below. If $\left\Vert u\right\Vert $ is large, then $\varrho(u)\geq2$, which yields
$\phi(\varrho(u))=0$, $\Psi(u)=\varrho(u)$. Thus $\Psi$ is coercive.

Let $\left\{  u_{n}\right\}  $ be a $\left(  PS\right)  _{c}$ sequence of
$\Psi$ with $c\leq0$, that is $\Psi(u_{n})\rightarrow c$, $\Psi^{\prime}%
(u_{n})\rightarrow0$. Then for $n$ large we have%
\begin{equation}
-\phi(\varrho(u_{n}))\int F(x,u_{n})=\Psi(u_{n})-\varrho(u_{n})\leq\frac{1}{2}%
-\varrho(u_{n})\text{.}\label{up}%
\end{equation}
We claim that%
\begin{equation}
1-\left(  \int F(x,u_{n})\right)  \phi^{\prime}(\varrho(u_{n}))\geq
1\text{.}\label{pu}%
\end{equation}
There are two cases:

\begin{enumerate}
\item If $\varrho(u_{n})<1$, then $\phi^{\prime}(\varrho(u_{n}))=0$ and (\ref{pu})
is true.

\item If $\varrho(u_{n})\geq1$, then (\ref{up}) and $\phi(\varrho(u_{n}))\geq0$ yield%
\[
-\phi(\varrho(u_{n}))\int F(x,u_{n})<0\text{,\qquad}\int F(x,u_{n})\geq0\text{.}%
\]
So we also have (\ref{pu}) becuase $\phi^{\prime}(\varrho(u_{n}))\leq0$.
\end{enumerate}

Since $\left\{  u_{n}\right\}  $ is bounded in $X$, we may assume that
$u_{n}\rightharpoonup u$ in $X$. Applying Proposition \ref{p1} we get%
\[
\int f(x,u_{n})(  u_{n}-u)  \rightarrow0\text{.}%
\]
This and $\langle\Psi^{\prime}(u_{n}),u_{n}-u\rangle\rightarrow0$, as well as (\ref{e0}) and
$|\phi(\varrho(u_{n}))|\le1$ imply%
\begin{align}
&  \hspace{-3em}\left(  1-\left(  \int F(x,u_{n})\right)  \phi^{\prime}%
(\varrho(u_{n}))\right)  \langle\varrho^{\prime}(u_{n}),u_{n}-u\rangle\nonumber\\
&  =\langle\Psi^{\prime}(u_{n}),u_{n}-u\rangle+\phi(\varrho(u_{n}))\int
f(x,u_{n})(  u_{n}-u)  \rightarrow0\text{.}%
\end{align}
Using (\ref{pu}) and $\varrho^{\prime}(u_{n})=\mathcal{L}u_{n}$, we get
\[
\langle\mathcal{L}u_{n},u_{n}-u\rangle\rightarrow0\text{.}%
\]
Because $\mathcal{L}$ is $\left(  S_{+}\right)  $-type, we deduce
$u_{n}\rightarrow u$ in $X$.\qed
\end{pf}

\paragraph{\it Proof of Theorem \ref{t2}} Because $\Psi(u)=\Phi(u)$ for $u\in\varrho^{-1}[0,1)$,
a neighborhood of $u=0$, it suffices to find a sequence of critical points
$\left\{  u_{n}\right\}  $ of $\Psi$ satisfying $\Psi(u_{n})\leq0$ and
$u_{n}\rightarrow 0$ in $X$ via Proposition \ref{pp1}. Having established
$\left(  PS\right)  _{c}$ for $c\leq0$ in Lemma \ref{l6}, it suffices to
verified condition (\ref{y}).

For each $k\in\mathbb{N}$ let $W_{k}$ be a $k$-dimensional subspace of $X$. There is $\Lambda_{k}%
\in\left(  0,\infty\right)  $ such that%
\begin{equation}
\left\vert u\right\vert _{p_{-}}^{p_{-}}\geq\Lambda_{k}\left\Vert u\right\Vert
^{p_{-}}\text{\qquad for }u\in W_{k}\text{.}\label{z}%
\end{equation}
By (\ref{f2}), there is $\delta>0$ such that%
\[
F(x,t)\geq\frac{2}{\Lambda_{k}p_{-}}\left\vert t\right\vert ^{p_{-}%
}\text{\qquad for }\left(  x,t\right)  \in\mathbb{R}^{N}\times\left(
-\delta,\delta\right)  \text{.}%
\]
Using $\dim W_k<\infty$ again, there is $\delta_{k}\in\left(  0,1\right)  $ such that if $u\in W_{k}$, $\left\Vert
u\right\Vert =\delta_{k}$, then $\left\vert u\right\vert _{\infty}<\delta$.

Now, if $u\in W_{k}\cap \partial B_{\delta_{k}}$, then for all $x\in\mathbb{R}^{N}$ we
have $\left\vert u(x)\right\vert <\delta$. Thus
\begin{align}
F(x,u(x))\geq\frac{2}{\Lambda_{k}p_{-}}\left\vert u(x)\right\vert ^{p_{-}%
}\text{.}\label{ed}
\end{align}
Noting that
\begin{equation}
\varrho(u)=\int\frac{1}{p(x)}\left(  \left\vert \nabla u\right\vert
^{p(x)}+V(x)\vert u\vert ^{p(x)}\right)  \leq\frac{1}{p_{-}%
}\left\Vert u\right\Vert ^{p_{-}}\text{,}\label{re}%
\end{equation}
using (\ref{z}), (\ref{ed}) and (\ref{re}) (which also implies $\phi(\varrho(u))=1$)), we get
\begin{align*}
\Psi(u)  & =\varrho(u)-\phi(\varrho(u))\int F(x,u)=\varrho(u)-\int F(x,u)\\
& \leq\frac{1}{p_{-}}\left\Vert u\right\Vert ^{p_{-}}-\frac{2}{\Lambda
_{k}p_{-}}\int\left\vert u\right\vert ^{p_{-}}=\frac{1}{p_{-}}\left\Vert
u\right\Vert ^{p_{-}}-\frac{2}{\Lambda_{k}p_{-}}\left\vert u\right\vert
_{p_{-}}^{p_{-}}\\
& \leq\left(  \frac{1}{p_{-}}-\frac{2}{p_{-}}\right)  \left\Vert u\right\Vert
^{p_{-}}=-\frac{\delta_{k}}{p_{-}}\text{.}%
\end{align*}
Consequently, condition (\ref{y}) is verified:
\[
\sup_{W_{k}\cap \partial B_{\rho_{k}}}\Psi\leq-\frac{\delta_{k}}{p_{-}}<0\text{.}\tag*{$\Box$}
\]

\begin{rem}
Replacing (\ref{f2}) in $(f_0)$ by the following much weaker condition:
\begin{align}
\lim_{\left\vert t\right\vert \rightarrow0}\frac{F(x,t)}{\left\vert
t\right\vert ^{p_{-}}}=+\infty\text{\quad point-wise in some ball }B_r(a)\subset\mathbb{R}%
^{N}\text{,} 
\end{align}
we can still verify (\ref{y}) as in the proof of \citet[Lemma 3.4]{liu2024multiple}. Therefore, Theorem \ref{t2} is still valid.
\end{rem}

\end{document}